\documentclass[12pt]{amsart}

\usepackage{amsmath,amssymb,amsfonts,amscd,enumerate,verbatim,calc,xypic}

\def\bB{{\boldsymbol B}}

\def\bF{{\boldsymbol F}}

\def\CB{{\mathcal B}}
\def\CC{{\mathcal C}}

\def\CF{{\mathcal F}}

\def\CK{{\mathcal K}}
\def\CM{{\mathcal M}}
\def\CP{{\mathcal P}}

\def\CT{{\mathcal T}}

\theoremstyle{plain}
\newtheorem{theorem}{Theorem}[section]
\newtheorem{lemma}[theorem]{Lemma}

\newtheorem{proposition}[theorem]{Proposition}

\newtheorem{conjecture}[theorem]{Conjecture}

\theoremstyle{definition}

\theoremstyle{remark}
\newtheorem{remark}[theorem]{Remark}
\newtheorem{chunk}[theorem]{}

\numberwithin{equation}{section}

\newcommand{\frgr}[1]{{\mathfrak g}_{Ad\rho_{}}}

\newcommand{\ds}[1]{\displaystyle{#1}}

\newcommand{\cal}[1]{\mathcal{#1}}
\newcommand{\braidp}[2]{P_{#1,#2}}
\newcommand{\braid}[2]{\overline{P}_{#1,#2}}

\newcommand{\incl}[1][r]{\hspace{0.2mm}\ar@<-0.22pc>@{^(-}[#1] \ar@<+0.22pc>@{-}[#1]}
\newcommand{\incld}[1][d]{\vspace{0.2mm}\ar@<-0.2pc>@{^(-}[#1] \ar@<+0.2pc>@{-}[#1]}

\def\CC{{\mathbb C}}
\def\FF{{\mathbb F}}

\def\PP{{\mathbb P}}

\def\ZZ{{\mathbb Z}}

\def\bar#1{\overline{#1}}

\def\into{\hookrightarrow}

\def\Mtwo{{\em Macaulay} 2\expandafter}

\begin{document}

\title[Mapping Class Groups]{Genus $2$ Mapping Class Groups are not K\"ahler}

\author{R\u{a}zvan Veliche}
\address{Department of Mathematics, Purdue University}
\curraddr{Department of Mathematics, University of Utah}
\email{rveliche@math.utah.edu}
\thanks{}


\date{}

\dedicatory{To Oana and ``AAA''}

\commby{}

\begin{abstract}
The goal of this note is to prove that the 
mapping class groups of closed orientable surfaces of genus 2 (with
punctures) are not K\"ahler. An application to compactifications of the moduli space of genus $g$ curves (with punctures) is given.
\end{abstract}

\maketitle

\section{Introduction}

It is a long-standing problem to determine which groups can occur in
the class $\CK$ of 
fundamental groups of compact K\"ahler manifolds. There are many
examples and positive results, but there are also many restrictions
(negative tests) that prohibit groups from being in $\CK$ (see
\cite{abc:abckt}).

Related to the class $\CK$, one can introduce several other classes of
fundamental groups: $\cal{QK}$ (of complements of divisors in K\"ahler
manifolds), $\CP$ (of projective varieties or manifolds) and $\cal{QP}$
(of quasi-projective varieties). Since every projective manifold is a
K\"ahler one, we have 

\begin{equation}
\xymatrix{
\xymatrixrowsep{1pc}
\stackrel{\ds \CP}{} \incl[r]\incl[d] & \stackrel{\ds \CK}{} \incl[d]\\
\cal{QP} \incl[r] & \cal{QK}\\
} 
\end{equation}
(and it is an open question whether the horizontal inclusions
are equalities).

The mapping class group $M(g,n)$ of (punctured) closed orientable surfaces is
in $\cal{QP}$ (\ref{mapgp_qkahler}). One may ask the question whether
this group is in $\CP$. In this note a negative answer is given in the
genus $2$ case.

An immediate consequence is that $\CM_{g,n}$, the moduli
space of $n$-punctured curves of genus $g$, cannot have a projective compactification with a ``thinner'' boundary (i.e. of higher codimension): 
the Satake compactification has a codimension two component ``at infinity'', and the
Deligne-Mumford compactification has a divisor there.

\section{Notations and Background}

\chunk
In what follows, a {\it variety} is an irreducible but possibly non-reduced scheme over $\CC$. A {\it manifold} is a smooth variety (i.e. it is reduced, too).
\chunk
For a compact orientable (connected) topological manifold $T$ one
defines the {\it   (pure) braid   group} 
=$\braidp{0}{n}(T)$ as the fundamental group of
$$\bF_{0,n}(T)=\{(z_1,\dots,z_n)\in T^n| 
z_i\neq z_j,\  \forall 1\leq i\neq j\leq n\}$$ There is a 
free action (by permutation) of the symmetric group $S_n$ on $\bF_{0,n}(T)$ and one
denotes the fundamental group of the quotient manifold
$\bB_{0,n}(T)=\bF_{0,n}(T)/S_n$ by $\braid{0}{n}(T)$
(called the {\it full braid group}).

\chunk
Introducing punctures, one defines
$\braidp{k}{n}(T)$ as $\pi_1(\bF_{0,n}(T\setminus Q_k))$ where $Q_k$ is
any set 
of distinct $k$ points in $T$ (its choice does not matter, as one can
always find a homeomorphism of $T$ mapping one choice onto another),
and respectively $\braid{k}{n}(T)$ as the fundamental group of
$\bB_{k,n}(T)=\bF_{k,n}(T)/S_n$.

\chunk
\label{PBfingen}
It is known that $P_{k,n}(T)$ and $B_{k,n}(T)$ are finitely generated groups; see e.g. \cite[\textsection 1.2,1.3]{birman:blm}.

\chunk
It is easy to see (\cite[Prop.1.1]{birman:blm}) that there are extensions: $$1\to \braidp{k}{n}(T)\to
\braid{k}{n}(T)\to S_n\to 1$$

\chunk
Related to these, one can consider $\CF_n(T)$
as the space of orientation-preserving homeomorphisms of $T$ fixing a
given $n$-tuple $(z^0_1,\dots,z^0_n)$ of distinct points of $T$, and
$\CB_n(T)$ as the space of orientation-preserving homeomorphisms
of $T$ fixing the set $\{z^0_1,\dots,z^0_n\}$. Both of these spaces
are endowed with the compact-open topology. Denote by
$M^p(T,n)=\pi_0(\CF_n(T))$ (the {\it pure mapping class group} of $T$)
and by
$M(T,n)=\pi_0(\CB_n(T))$ the {\it (full) mapping class group} of
$T$ (also denoted in the literature by $\Gamma_{g,n}$). These are the
groups of isotopy classes of (pointwise, respectively setwise) puncture-preserving diffeomorphisms.

\chunk
\label{MpinM}
As in the case of the braid groups, we have the extensions: 

$$1\to M^p(T,n)\to M(T,n)\to S_n\to 1$$

\chunk
In the case when $T$ is a closed oriented surface of genus $g$,
$M(T,n)$ is classically denoted by $M(g,n)$ (\cite{birman:blm}) or $M_g^n$.

\chunk
This note is concerned with the structure of $M(g,n)$, and
particularly with the possibility of it being the fundamental group of
a K\"ahler manifold. 

\chunk
\label{mapgp_qkahler}
Note that $M(g,n)$ is the fundamental group
of a quasi-projective manifold. Namely, the Teichm\"uller space $\CT_{g,n}$ of 
equivalence classes of orientation (and puncture) preserving
diffeomorphisms of a genus $g$ smooth curve, is a contractible complex
manifold. The 
mapping class group $M(g,n)$ acts properly discontinuously on $\CT_{g,n}$, but the action is
not necessarily free; actually, even the pure mapping class (sub)group
$M^p(g,n)$ might act on
$\CT_{g,n}$ with fixed points (this does not happen as soon as $n>2g+2$).  However, the
quotient $\CM$ has at most finite-quotient singularities (it is an
orbifold); $M(g,n)$ becomes its orbifold fundamental group $\pi_1^{orb}(\CM)$,
and there exists a finite index subgroup $H$  
(which can be chosen to be normal) of $M(g,n)$ which acts freely on
$\CT_{g,n}$ (\cite[12.8]{hain:moduli}), and such that $H\backslash \CT_{g,n}$ is quasi-projective and smooth. That means that there is a finite (orbifold) cover
$\widetilde{\CM}\to \CM$ with $\widetilde{\CM}$ a smooth quasi-projective manifold.
Now one may use the argument in \cite[Lemma 1.15, p.7]{abc:abckt} to
show that $M(g,n)$ is the fundamental group of a (smooth) quasi-projective manifold.

\chunk
We recall the definition of commensurability: it is the equivalence
relation in the category of groups generated by morphisms with finite
kernel and cokernel.


\chunk
\label{ends}
We will be using some results on the ends of a (K\"ahler) group, so
it's worth saying a few words on the subject.

For a topological space $X$, the set of ends of $X$ is the inverse
limit $\lim \pi_0(X\setminus K)$, where the limit is taken as $K$
ranges over all compact subsets of $X$. The cardinality of this set is
denoted by $e(X)$. If a group $G$ acts freely on a connected
simplicial complex $X$ with finite quotient, then $e(X)$ is
independent of $X$, and is denoted by $e(G)$; this is the number of
ends of the group $G$.

\chunk
\label{facts_ends}
Facts (see \cite{cohen:grcd1}, \cite{scott:top.met.gr}):
A free group with more than two generators has infinitely many
ends. The number of ends of an infinite group can be 1,2 or
$\infty$. Commensurable groups have the same number of ends.

\chunk
\label{facts_kahler}
A finite index subgroup of a group in $\CK$ (or $\CP$) is also in $\CK$ (resp. $\CP$); this is
roughly due to the possibility of lifting and averaging a K\"ahler metric to a
finite unramified cover, respectively to the fact that the pullback of an ample divisor to a finite cover is ample.

However, a finite extension of a group in $\CK$ need not be in $\CK$; e.g. if $G\cong \ZZ^2\rtimes \ZZ_2$ by the action switching the generators of $\ZZ^2$, then $G$ is not in $\CK$ (the first Betti number is $1$).

\section{Preliminaries}

We start by a minor generalization of a result in 
\cite[(4.1, p.151)]{birman:blm}:

\begin{lemma}
The evaluation map $\epsilon_{g,n,k}:\CF_k(T)\to F_{k,n}(T)$ is
a locally trivial fibering with fibre $\CF_{n+k}(T)$.
\end{lemma}

\begin{proof}
The argument is quasi-identical to the proof in \cite{birman:blm}, but we
include it here for the reader's convenience. 

Fix an $(n+k)$-tuple $(z^0_1,\dots,z^0_{k+n})$ of distinct points in
$T$. Without restricting the generality, the homeomorphisms in
$\CF_i(M)$ can be taken to fix $(z^0_1,\dots,z^0_i)$ (for $i\in
\{k,n+k\}$). $\CF_{n+k}(T)$ becomes then a closed
topological subgroup of $\CF_k(T)$.

The
evaluation map $\epsilon_{g,n,k}$ takes a homeomorphism $h$ to 
$(h(z^0_{k+1}),\dots,h(z^0_{n+k}))$, which is clearly a point in $F_{k,n}(T)$ thought of as
$F_{0,n}(T\setminus\{z^0_1,\dots,z^0_k\})$.

$\epsilon_{g,k,n}$ is surjective (any finite set of points can be
taken to another set of the same cardinality by a homeomorphism of
$T$) and two homeomorphisms $h_1$ and $h_2$ have the same image under
$\epsilon_{g,k,n}$ if and only if they are in the same left coset of
$\CF_{n+k}(T)$ in $\CF_k(T)$.

Thus $F_{k,n}(T)$ is naturally (topologically) identified with the
quotient space \\
$\CF_k(T)/\CF_{n+k}(T)$, and $\epsilon_{g,k,n}$ becomes the projection
map which turns $F_{k,n}(T)$ into a homogeneous space.

By \cite[p.99]{hu:ho}, the theorem follows as soon as we show the
existence of a local 
section for $\epsilon_{g,k,n}$ near ${\bf
  z}^0=\epsilon_{g,k,n}(\CF_{n+k}(T))$.  

For that, choose pairwise
disjoint Euclidean neighbourhoods $U(z^0_{k+1}),\dots,U(z^0_{n+k})$ of
\\ $z^0_{k+1},\dots,z^0_{k+n}$; we require that the union of these neighbourhoods
doesn't intersect $\{z^0_1,\dots, z^0_k\}$. 

Then $U({\bf
  z}^0)=\{(u_{k+1},\dots,u_{k+n})\,|\, u_i\in U(z^0_i)\}$ is a
neighbourhood of ${\bf z}^0$ in $F_{k,n}(T)$. Construct a family of
homeomorphisms $\{{\it f}_u\in \CF_k(T)\,|\, u\in U({\bf z}^0)\}$
depending continuously on $u$ and such that for each $u\in U({\bf
  u}^0)$:

\begin{itemize}
\item ${\it f}_u(z^0_{k+j})=u_{k+j}$ ($j=\overline{1,n}$) and
\item ${\it  f}_u\left(T\setminus 
\bigcup\limits_{j=k+1}^{k+n} U(z^0_j)\right)$
is the identity. 
\end{itemize}

Such a family exists because, for example, there exists a continuous
family of homeomorphisms (even diffeomorphisms) deforming the center
of a ball to any of the points in the ball.

Now the section we wanted is given by $s(u)={\it f}_u$.
\end{proof}

\chunk
\label{Myseq}
In the case of a closed oriented surface of genus $g$, the lemma
gives, via the long exact sequence in homotopy associated to 
a fibration:

$$P_{k,n}(T_g)\to M^p(g,n+k)\to M^p(g,k)\to \pi_0(F_{k,n}(T))=1$$

\chunk
\label{M2toM0}
We will make essential use of the fact that there is an exact
sequence:

$$1\to \ZZ_2\to M(2,0)\to M(0,6)\to 1 $$

(see \cite[p.188]{birman:blm}).

\chunk
\label{Mseq}
In \cite[Thm. 4.2, p.152]{birman:blm} it is proven that there is an exact
sequence:

$$1\to P_{0,n}(T_g)\to M^p(g,n)\to M^p(g,0)\to 1$$
for $g\geq 2$, and a similar one:

$$1\to P_{0,n}(T_g)/\mbox{center}\to M^p(g,n)\to M^p(g,0)\to 1$$
if $g=1$, $n\geq 2$ or $g=0$, $n\geq 3$.

\chunk
\label{Pseq}
In \cite[Thm 1.4, p.14]{birman:blm} the existence of the following exact
sequence is proven:

$$1\to P_{n-1,1}(M)\to P_{0,n}(M)\to P_{0,n-1}(M)\to 1$$

\noindent
for any $n\geq 1$.
One needs $M$ to satisfy $\pi_i(M\setminus Q_m)=1$ for $i=0,2,3$ and
for every $m\geq 0$ as hypotheses. 
In the case of $S^2$ (which fails the above test for $m=0$) the exact
sequence is valid only for $n\geq 4$ (see \cite[p.34]{birman:blm}).

\chunk
\label{aboutS2}
We will need the fact that $P_{0,3}(S^2)=\ZZ_2$; see
\cite[p.34]{birman:blm}.

\section{Main result}

\begin{theorem}
\label{main_thm}
The groups $M(2,n)$ are not K\"ahler, for any $n\geq 0$.
\end{theorem}

\begin{proof}
The idea of the proof is rather simple: show that these groups have
too ``big'' an image to be K\"ahler.

More precisely, \cite{arapura:ABR} prove that an extension of a group with
infinitely many ends by a finitely generated group cannot be K\"ahler.

Thus, if we prove that inside each of the $M(2,n)$'s there is a
finite index subgroup mapping (with finitely generated kernel) onto a
group with infinitely many ends, we're done.

The natural subgroup of finite index to consider in each $M(g,n)$ is
$M^p(g,n)$ (\ref{MpinM}). The 
exact sequence (\ref{Mseq}) for $g=2$ and the remark (\ref{PBfingen}) show
that it is enough to prove that $M^p(2,0)=M(2,0)$ maps (with finitely
generated kernel) onto a group with infinitely many ends.

Now (\ref{M2toM0}) reduces the problem to $M(0,6)$. This doesn't map
directly to what we need, but it's (finite index) subgroup $M^p(0,6)$ maps, by
(\ref{Myseq}), onto $M^p(0,4)$ with finitely generated kernel.

Now the second sequence in (\ref{Mseq}) with $g=0$ and $n=4$ gives
$M^p(0,4)\cong P_{0,4}(S^2)/\mbox{center}$.

But (\ref{Pseq}) for $n=4$ (and $g=0$) gives:

$$1\to P_{3,1}(S^2)\to P_{0,4}(S^2)\to P_{0,3}(S^2)\to 1$$

But $P_{0,3}(S^2)=\ZZ_2$ (\ref{aboutS2}), and $P_{3,1}(S^2)=\pi_1(S^2\setminus
\{z^0_1,z^0_2,z^0_3\})\cong\pi_1(\CC\setminus{0,1})=\FF_2$ (the free
group on two generators). Now it is easy to see that the center of
$P_{0,4}(S^2)$ is either trivial or isomorphic to $\ZZ_2$ (it is
actually the latter, see \cite[Lemma 4.2.3, p.154]{birman:blm}), and
since it cannot intersect the $\FF_2$, one has a finite index copy of
$\FF_2$ inside $M^p(0,4)$ (we actually have $M^p(0,4)\cong \FF_2$). 
Pulling back this copy of $\FF_2$ to either $M^p(0,6)$ or $M^p(2,0)$
finishes the proof.

\end{proof}

\section{Higher genus; Conjecture}

\chunk
The exact sequence 
(\ref{Myseq}) for compact surfaces of genus $g\geq 3$, $n=1$ and
$k\geq 1$ yields:

$$P_{k,1}(T_g)\stackrel{d}{\to} M^p(g,k+1)\stackrel{i}{\to}
M^p(g,k)\to 1$$ 

\noindent
where $d$ is the ``boundary'' map $\pi_1(F_{1,1}(T_g)\to
\pi_0(\CF_2(T_g))$ and $i$ is induced by the fibre-inclusion
$\CF_2(T_g)\to \CF_1(T_g)$.

\begin{lemma} $\ker(d)\subseteq \mbox{center}( P_{k,1}(T_g))$
\end{lemma}
\begin{proof} Easy adaptation of the proof of \cite[Lemma 4.2.1,
  p.153]{birman:blm}. 
\end{proof}

Since $P_{k,1}(T_g)=\pi_1(T_g\setminus\{z_1,\dots,z_k\})=\FF_{2g+k-1}$,
the lemma
 gives the existence of a finitely generated free normal subgroup
of infinite index in $M^p(g,k+1)$.

\chunk
At this point, we would like to introduce the class ${\bf NNF}$ of
groups that do {\bf not} have {\bf normal free finitely generated subgroups} (of
rank at least $2$).

\begin{conjecture}
\label{conj_kahler}
 $\CK\subset {\bf NNF}$
\end{conjecture}

\begin{remark}
The validity of this conjecture would imply that all
the mapping class groups $M(g,n)$ with $g\geq 2$ and $n\geq 2$ are not
in $\CK$. 
\end{remark}

\begin{remark}
\label{finind}
It is known that any (finitely generated) free group has a finite
index normal subgroup of odd rank. Then (\ref{facts_kahler}) shows
that groups in $\CK$ cannot have finite index normal free
subgroups (these would be finitely generated, actually).

Thus the conjecture refers to the infinite index case.
\end{remark}

\chunk
As support for this conjecture we cite the work of
\cite{catanese:ninf}, showing that the fundamental groups of surfaces
of genus $g$ are of class ${\bf NNI}$. That is a more restrictive
class, introduced by Catanese, where {\bf no finitely generated normal
subgroups} (not necessarily free) of {\bf infinite index}
are allowed in the given group.

${\bf NNI}\nsubseteq {\bf NNF}$ (take for example a direct
product of a free group, of rank at least $2$, with $\ZZ_2$). But if
one enlarges {\bf NNF} to {\bf NNIF} (prohibiting only {\bf infinite
  index} finitely generated normal free subgroups, of rank at least $2$), one has ${\bf
  NNI}\subset {\bf NNIF}$.

 Note that we have $\cal{K}\cap {\bf NNI}\subseteq \CK\cap {\bf
   NNIF}=\CK\cap {\bf 
   NNF}$ (the latter by (\ref{facts_kahler})). The conjecture simply
 says  that $\CK\subset {\bf NNIF}$.

\begin{remark} The classes {\bf NNI}, {\bf NNIF}, $\cal{K}\cap {\bf
    NNI}$ and $\CK\cap {\bf NNIF}$ are not closed under direct
  products. It is an easy observation 
  that {\bf NNF} and $\CK\cap {\bf NNF}$ are closed under this
  operation. This implies that  arbitrary (finite) direct products of
  finite groups and fundamental groups of smooth projective curves are in
  $\CK \cap {\bf NNF}$.
\end{remark}

\begin{remark}
As proved in \cite[Prop.2.5,
p.23]{cohen:grcd1}, if a group $G$ has a subgroup $H$ which is not
locally finite, but which is included in a finitely generated
subgroup of infinite index, then $G$ has one end. In particular a
direct product of two free groups of rank greater than $2$ has one end
only. So the obstruction in
(and a proof of) the above conjecture should have more geometric content.
\end{remark}

\section{Relevance to Compactifications of $\CM_{g,n}$}

\chunk 
A natural question that arises in connnection to $\CM_{g,n}$ is
whether one could find a projective (more 
generally, a weakly K\"ahlerian) compactification for it, such that
the ``boundary'' has high codimension. The Deligne-Mumford
compactification has a divisor, and the Satake
compactification has a codimension two component in the boundary.

\begin{proposition}$\CM_{2,n}$ has no projective compactification with
  boundary of codimension at least $3$.
\end{proposition}

Using \ref{main_thm}, \ref{mapgp_qkahler} and \ref{facts_kahler}, this is an immediate consequence of the following:

\begin{lemma}
Given:
\begin{enumerate}[\rm\quad (1)]
\item a quasi-projective orbifold $X$ with an irreducible universal cover $\widetilde{X}$
\item a (normal) finite index subgroup $H$ of  $\pi_1^{orb}(X)$ such that $Y=H\backslash \widetilde{X}$ is quasi-projective and smooth, and such that 
\item $H\notin \CP$
\end{enumerate}
then $X$ cannot have a projective compactification with boundary of codimension at
least $3$.

\end{lemma}

\begin{proof}

Assume there would be such a compactification, denoted by
$\bar{X}\into \PP^N$.
The map $p: Y\to X$ is onto and finite by hypothesis. Composition with the
inclusion $X\into \bar{X}$ gives a dominant quasi-finite map
$f: Y\to \bar{X}\subseteq \PP^N$. 

Since the codimension of the boundary $\bar{X}\setminus X$ is at least $3$ and $\bar{X}$ is projective, we can find an irreducible
subscheme $L'=L\cap \bar{X}$ in $\bar{X}$, of dimension $2$, with
$L$ a linear space of $\PP^N$, and such that $L'\cap
(\bar{X}\setminus X)=\emptyset$. This implies that for a sufficiently
small $\epsilon$ (and a Riemannian metric on $\PP^N$), the neighbourhood 
$L'_\epsilon$ of points at distance less than $\epsilon$ from $L$ is also
avoiding the boundary $\bar{X}\setminus X$.

Now we can use
the theorem of Goresky and MacPherson (\cite{goresky:strat}, see also
\cite[\S 9]{fulton:connect}) regarding pullbacks of linear subspaces
of $\PP^N$ via quasi-finite morphisms, to get that
$\pi_i(Y,f^{-1}(L'_\epsilon))=1$ for $i\leq 2$.

This implies
$\pi_1(Y)=\pi_1(f^{-1}(L'_\epsilon))=\pi_1(f^{-1}(L'))$.  

Now $f^{-1}(L')=p^{-1}(L')\to L'$ is a finite map, and $L'$ is projective.  Then $p^{-1}(L')$ is also projective. 

But then $\pi_1(Y)=H$ would be in $\CP$, contradiction.

\end{proof}

Note: it is sufficient to assume in the lemma that $Y$ is a (connected, but possibly reduced) local complete intersection (see \cite[\S 9]{fulton:connect}).

\bibliographystyle{amsplain}

\end{document}